\documentclass[letterpaper,aps,twocolumn,floatfix]{article}
\usepackage{yfonts}
\usepackage{enumitem}
\usepackage{graphicx}
\usepackage{epsfig}
\usepackage{amsmath}
\usepackage{amssymb}
\usepackage{amstext}
\usepackage{amsfonts}
\usepackage{hyperref}
\usepackage{alltt}
\usepackage{setspace}
\usepackage{euscript}
\usepackage{yfonts}
\usepackage{accents}
\usepackage{psfrag}
\usepackage{array}
\usepackage{abstract}
\newcommand{\old}[1]{{}} \old{}

\begin{document}

\twocolumn[

\title{Non-Smooth Solutions of the Navier-Stokes Equation and their Means}

\author{J. Glimm \\
Stony Brook University, Stony Brook, NY 11794,
and GlimmAnalytics LLC, USA  \\
glimm@ams.sunysb.edu
\and
J. Petrillo \\
Point72, 55 Hudson Yards, New York, NY 10001 \\
Jarret.petrillo@point72.com }

\maketitle

\let\oldclearpage\clearpage
\renewcommand{\clearpage}{}

\begin{@twocolumnfalse}

\begin{abstract}

Non-smooth (finite time blowup) Leray-Hopf solutions of the incompressible
Navier-Stokes equation are constructed.
The initial data for blowup is characterized by  nonzero energy related
turbulent fluctuations.
The construction occurs in a finite periodic cube $\mathbb{T}^3$.

The mean value of a weak solution of the Navier-Stokes equation is identified
as a smooth solution of the Navier-Stokes equation. 
\end{abstract}
\end{@twocolumnfalse}

\let\clearpage\oldclearpage

]

\tableofcontents

\section{Introduction}
\label{sec:intro}

\subsection{Main Results}
\label{sec:main}

The main result of this paper is the
construction of non-smooth solutions of the Navier-Stokes equation
of the Leray-Hopf (LH) class of weak solutions. The initial data giving
rise to the non-smooth solutions is turbulent. Turbulent initial data
is identified as nonzero turbulent fluctuations. Within the limits of
the present analysis, these fluctuations are energy (not enstrophy) related.

Initial data which is not turbulent is characterized as the mean value of a
general Navier-Stokes solution. These solutions are entropy minimizing.

All proofs are made easier by working in $\mathcal{V}^*$ rather than in
$\mathcal{H}$.

The achievability of the initial data from limits
of the small time data follows from known results.

Analyticity in $\mathcal{V}^*$ of the weak solutions of the Navier-Stokes
equation is elementary. 
Blowup in finite time for selected turbulent initial conditions is established. 

A concluding section characterizes the mean of a weak solution of the
Navier-Stokes equation as a smooth solution of the Navier-Stokes equation.

In this paper, only incompressible isothermal single fluids (not mixtures)
are considered.  The fluids are real (not complex) valued. They
are defined in a periodic spatial domain $\mathbb{T}^3$.

\subsection{The Role of the LH Assumption}
\label{sec:role-lh}

The LH assumption has a unique role in this paper, which is to guarantee
existence of weak solutions whose initial conditions have been specified.
This role occurs uniquely in the Blowup section, Sec.~\ref{sec:blowup}. 

The LH assumption is not even needed for the proof that the solution
has intervals of regularity separated by possible singular times $t_j$.
This property can be deduced from the analyticity property, proven
for general weak solutions in Sec.~\ref{sec:III}.

\subsection{Prior Results}
\label{sec:prior}

A review of Navier-Stokes solutions can be found in \cite{BedVic22}. 
The criteria of \cite{BeaKatMaj84} are necessary for blowup. The numerical
study \cite{Hou23} with additional numerical and theoretical
references cited provide strong evidence for the
existence of non-smooth solutions of the Navier-Stokes equation.

Weak solutions for the Euler equation
have been constructed, with weakly defined limits and function spaces
\cite{DiPMaj87,BreDeLSzk11,BreFei18,GanWes09,BucVic19,KukVicWan20}
and additional references cited in these sources.
See \cite{Ber23} and references cited there for
a discussion of the Navier-Stokes and Euler equations
and their relation to the Onsager conjecture.

The Foias \cite{Foi74,FoiManRos04}
description of Navier-Stokes turbulence as an LH weak solution
provides the constructive step needed in our analysis. Other ideas, including
analyticity, from this body of work are also used.
The Robinson analysis  \cite{RobRodSad16}, Ch. 3, is more general and applies
to arbitrary weak solutions of the Navier-Stokes equation. 

The existence of entropy production maximizing solutions of the
Navier-Stokes equation is established in \cite{GliLazChe20a,CheGliSai24b}.
This result is not sufficient for our needs and will be extended.

The theory of vector and tensor Fourier analysis 
\cite{VarMosKhe88,IvePhi08} is used.

The Lagrangian for space time smooth fluids, derived from
\cite{Jon15,FukFuj10,FukFuj12,HolMarRat98,Tay11,SelWhi68},
has influenced our thinking.
Scaling analysis, based on the renormalization group has also influenced
our thinking. Of the many papers on this subject,
we mention \cite{BogShi59,MonYag75,Wil83}

Additional references are added as needed throughout the text.
Additional references can be found within the sources cited above.
This paper is not a scholarly review.

\subsection{The Fluid Millennium Prize}
\label{sec:mill}

A technical wording of this specification is given by
C. Fefferman in several options. The spatially periodic option B is:

Take $\nu > 0$ and $n = 3$. Let $u^0(x)$ be any smooth, divergence-free
vector field satisfying
\begin{itemize}
\item smoothness, divergence free, periodicity of $u^0$.
\item No restrictions are placed on $f$ other than periodicity.
\end{itemize}
Then there exist smooth functions $p(x,t)$,
$u_i(x,t)$ on $\mathbb{R}^3 \times [0,\infty)$
that satisfy conditions including

1. the Navier-Stokes equation;
2. divergence free;
3. initial conditions given by $u^0(x)$;
10. periodicity;
11. smoothness.

The result, if achieved, is a solution of the prize problem in the 
positive. If it can be shown that these conditions are not satisfied for
some weak solution of the Navier-Stokes equation, the prize problem is
solved in the negative.

We solve the prize problem in the negative.
We choose $f = 0$.
Initial conditions $u^0$ that are turbulent in the energy analysis
space allow blow up in finite time.
%The prize problem is solved in the negative.

This paper is offered as a solution of the Millennium Fluids problem.

We show that the mean of any weak solution is smooth for all times.
The initial data for the mean is the mean of the weak solution. 
The mean does not have general initial data.
For this reason, the mean is not a positive solution of the Millennium Fluids
problem.

\subsection{The Reason for Blowup}
\label{sec:reason}

\paragraph*{Notation for Weak Limits}

$L^p(0,T;\mathcal{V}^*)$ does not describe weak convergence
by test functions $\phi \in L^1(0,T;\mathcal{V})$. Weak convergence
occurs in an $L^p$ space
\begin{equation}
\label{eq:weak-conv}
\mathcal{W}^{p} = L^{p}(0,T;\mathcal{V}^*L^1(0,T;\mathcal{V})) \ .
\end{equation}

The space $L^1(0,T;\mathcal{V})$ which multiplies $\mathcal{V}^*$
is a space of cutoff functions $\phi$ needed for the weak convergence
to occur. See Sec.~\ref{sec:cutoff}.

Where the context is clear, the Lebesgue time exponent $p$ will be omitted.

The Navier-Stokes solutions $u \in \mathcal{W}^p$ decay with the rate
$O(e^{-\nu t/p})$ or faster. When speaking of decay rates, it is 
as a minimum, or lower bound to the actual decay, with an important
exception: the SRI $O(e^{-\nu t/2})$ decay rate is exact.
It is not an upper bound.

\paragraph*{Three Contradictory Facts}

Assume $u$ is entropy production maximizing.

Three facts are logicly contradictory.

The contradiction is resolved by blowup and a loss of solution smoothness.

\begin{enumerate}
\item $u$ is analytic in $\mathcal{W}^2$ for all $t > 0$. 
With SRI, $u$ has an exact slow rate dissipation
(not an upper bound dissipation rate)
$O(e^{-\nu t/2})$ as an element of $\mathcal{W}^2$.

{\it Proof}: See Lemma~\ref{sec:ent}.1. 

\item The solution fluctuation $\nu_t$ defined by (\ref{eq:nu-NL-int}) is
analytic in $\mathcal{W}^{4/3}$ for all $t > 0$
for entropy production maximizing (and dissipation rate maximizing)
solutions.
$\nu_t$ decays at the rapid rate $O(e^{-3\nu t/4})$ as an element of  
$\mathcal{W}^{4/3}$.

{\it Proof}: See Theorem~\ref{sec:ent}.1 (b).

\item For entropy production maximizing (and dissipation rate maximizing)
solutions, the rapidly decreasing $\nu_t$ value in $\mathcal{W}^{4/3}$
is larger than the slowly decreasing $u$ value in $\mathcal{W}^2$.
That is, 
$\mathcal{W}^{4/3} \supseteq \mathcal{W}^2$.

{\it Proof}: See Theorem~\ref{sec:ent}.1 (c).

\end{enumerate}

\begin{figure}
\center{\includegraphics[width=225pt]{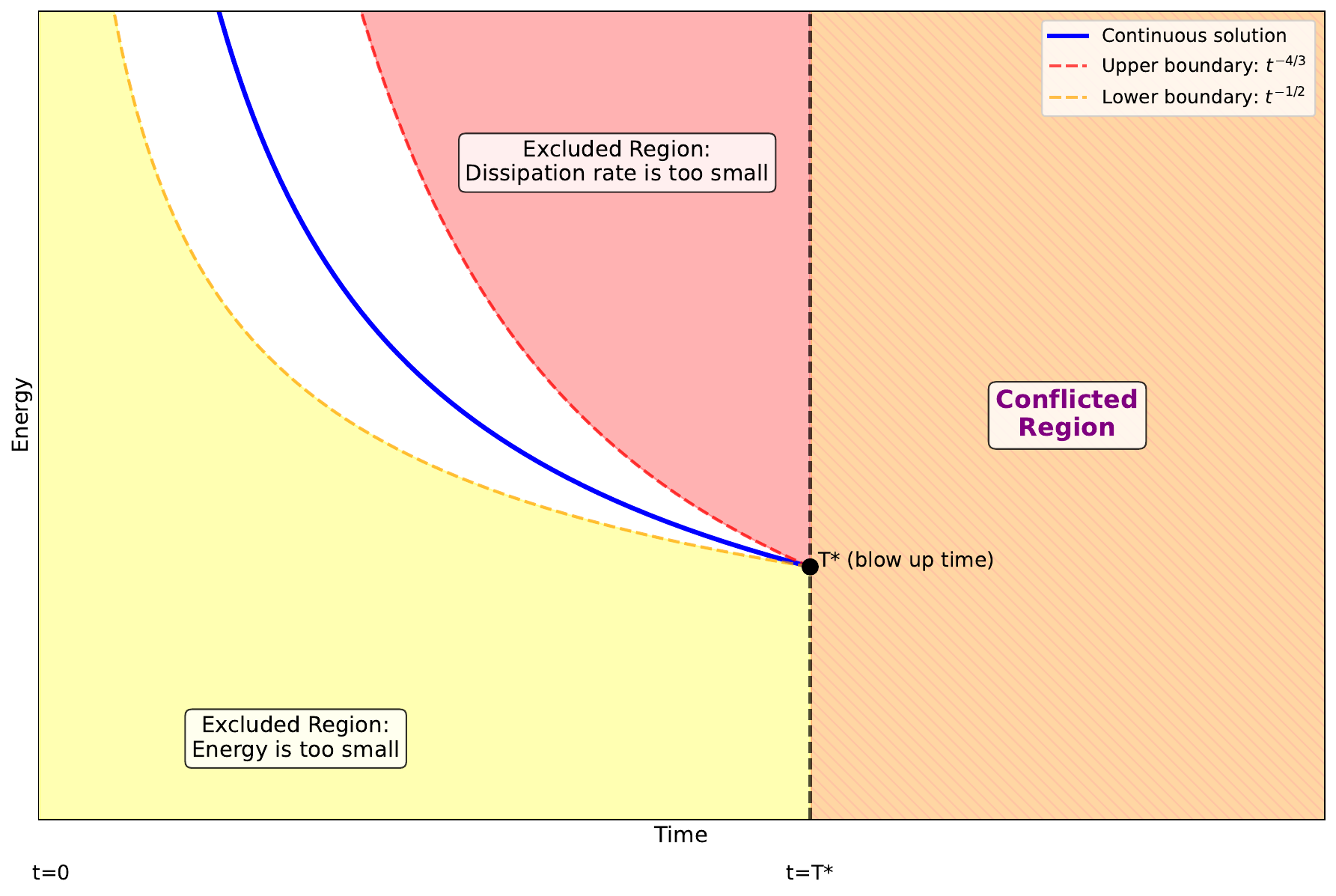}}
\caption{
\label{fig:blowup}
Two decay laws drive the blowup, one for the solution viscous decay 
$O(t^{-1/2})$ and the other for the decay $O(t^{-4/3})$ of fluctuations.
Each introduces a forbidden region, with a bound independent of time $t > 0$.
The forbidden solution decay region is shown in yellow
across the bottom of the figure. The turbulent fluctuation decay region is 
the upper region, shown in red. It is required that the solution decay
curve is below the fluctuation decay curve.
The two excluded regions overlap at a time $T^*$ in an excluded orange region.
The smooth solution cannot enter the excluded region. At
$t = T^*$, the solution must violate its fluctuation decay curve
bound, and by a discontinuous
amount. It is not smooth, and at the finite time $T^*$, it blows up.
}
\end{figure}

This impossibility is illustrated in Fig.~\ref{fig:blowup}.

The contradiction expressed in terms of decay rates,
using the exact dissipation rate for the SRI $u$ decay
and the lower bound for the $p = 4/3$ dissipation,
\begin{equation}
\label{eq:compare}
O(e^{-3\nu T/4}) \geq O(e^{-\nu T/2}) \ ,
\end{equation}
is logicly contradictory for large $T$.

\section{Preliminaries}
\label{sec:tech}

\subsection{Function Spaces}
\label{sec:fun}

\paragraph*{Configuration Space}

The Hilbert space for this theory is the space $\mathcal{H}$ of
$L_2$ incompressible  vector fields defined on the periodic cube $\mathbb{T}^3$,
also called configuration space. The inner product in this space is denoted
$\langle \cdot, \cdot \rangle$, and referred to as the bare vacuum inner 
product. The associated $L_1$ space functional is denoted 
$\langle \cdot \rangle$ (with no comma).

\paragraph*{Powers of the Stokes Operator}

The Stokes operator is denoted $A$.  For space periodic boundary conditions,
\begin{equation}
\label{eq:stokes}
A = -\nu \Delta \ ,
\end{equation}
where $\Delta$ is the Laplace operator, \cite{FoiManRos04}, Ch. II.6, Eq. (6.16).

The fractional power of the Stokes operator is $-\nu \Delta^s $.
Using these powers, we define fractional Hilbert spaces 
$\mathcal{H}_s$ with the inner product 
$\langle u,v\rangle_{\mathcal{H}_s} = \langle (u,(-\nu \Delta^s) v\rangle$.
Incompressibility is imposed as part of the definition of $\mathcal{H}_s$.

\paragraph*{The Leray Projector}

Incompressibility, $\sum_{i = 1}^3 \partial_ i u_i = 0$, is imposed on
all vector fields.

Define the Leray projector $\mathbb{P}$ as the projection onto the
divergence free subspace of velocity fields, acting on an $L_2$ space. 

$\mathbb{P}_q$ denotes the Leray projector
acting on an $L_q$ space, with an adjoint 
$\mathbb{P}_p$ acting on the dual $L_p$ space.
See \cite{Eva10}.

{\bf The Spaces $\mathcal{V}$ and $\mathcal{V}'$}

The space $\mathcal{H}_1$ is denoted $\mathcal{V}.$ Its dual is
$\mathcal{V}' = \mathcal{H}_{-1}$.

\subsection{Euler and Navier-Stokes Equations}
\label{sec:ns}

The Euler and Navier-Stokes equations
are defined physically as the law of conservation of momentum.
The weak solutions of these equations are interpreted through space time
integration, meaning that
they are valid as an expectation value with a cutoff function $\phi$,
as in (\ref{eq:energy-weak}).

The existence of nonunique and nonphysical
solutions of these equations have been known for the 1D Euler equation
since the days of Riemann, for
the 2D Euler equation \cite{Sch93} and in 3D even for the Navier-Stokes
equation; see \cite{BucVic19,BucVic20} and references cited there.

Evaluation of a divergent integral in a Sobolev space of negative index
does not change the theory, but modifies the
function spaces in which it is described. In view of the
infinities of turbulence, these negative index Sobolev
spaces are essential and can never be removed. In contrast, a cutoff, if used,
is temporary and removed in the sense of limits.

\paragraph*{Convective Derivatives}

The convective derivative is defined as
\begin{equation}
\label{eq:DDT}
\frac{D}{Dt} =  \frac{\partial}{\partial t} + u \cdot \nabla  \ .
\end{equation}

The Euler equation, expressed using Lagrangian derivatives, is
\begin{equation}
\label{eq:L-E}
\mathbb{P} \frac{Du_i}{Dt} = f_i  \ ,
\end{equation}
with $f$ the stirring force, assumed divergence-free.

\paragraph*{ The Stirring Force f}

Restrictions on $f$ are based on physics, as unrestricted stirring is
consistent mathematically with the Euler and Navier-Stokes equations.
The choice $f = 0$ is made within this paper.

\paragraph*{Stress Rates}

$S_{ij}$ is the rate of stress tensor, defined as
\begin{equation}
\label{eq:S}
S_{ij} = \frac{\partial u_i}{\partial x_j} \ .
\end{equation}
Vorticity $\omega$ is defined as
\begin{equation}
\label{eq:vor}
\omega_i = (\nabla \times u)_i 
\end{equation}
and the vorticity related rate of stress is
\begin{equation}
\label{eq:vor-S}
\nabla \times S_{ij} = \frac{\partial \omega_i}{\partial x_j} \ .
\end{equation}

The stress-strain tensor
\begin{equation}
\label{eq:sigma}
\sigma_{ij} = \frac{1}{2} \Big ( S_{ij} +  S_{ji} \Big )
\end{equation}
specifies the fluids modeled, an incompressible isotropic Newtonian fluid.

The Navier-Stokes equation written using Lagrangian derivatives is
\begin{equation}
\label{eq:L-NS}
\mathbb{P} \frac{Du_i}{Dt} = \nu \mathbb{P}( \partial_j \cdot S_{ij}) + f_i  \ .
\end{equation}

The Eulerian Navier-Stokes equation,
\begin{equation}
\label{eq:E-NS}
\bigg (\frac{\partial u_i}{\partial t} +
\mathbb{P} \Big(u_j \frac{\partial u_i}{\partial x_j} \Big) \bigg)
=\nu \mathbb{P} (\nabla_j \cdot \sigma_{ij}(x,t)) +f_i \ ,
\end{equation}
is formally equivalent to (\ref{eq:L-NS}) with
$\nu = 0$ to define the Euler equation.

\subsection{Cutoffs}
\label{sec:cutoff}

Cutoffs are defined as a product of a spatial and a temporal cutoff.
The spatial cutoff is the projection of the velocity field $u$ onto a finite
dimensional space of smooth basis elements of $\mathcal{H}$.
The temporal cutoff is defined by convolution by a smooth test function,
which is time retarded in the sense of being supported in $t < 0$.

A cutoff is defined as a modification to the theory.
The cutoff theory
restricts the spatial integration to a finite dimensional subspace of the
configuration space $\mathcal{H}$.

The complete criteria to define a weak solution of the Navier-Stokes equation
are formulated in terms of cutoffs.
Cutoffs are defined as functions $\phi$ belonging to 
\begin{equation}
\label{eq:d}
\begin{split}
& \widetilde{\mathcal{D}}_\sigma (\mathbb{T}^3) 
 = \Big \{ \phi \; : \; \phi = \sum_{k = 1}^N \alpha_k(t) a_k(x), \\
 & \mathrm{Div} \phi = 0,  \\
& \alpha_k \in C_c^1([0,\infty )), \; a_k \in \mathcal{N}, \; N \in \mathbb{N} 
\Big \}  \ .\\
\end{split}
\end{equation}
Here $\mathcal{C}^1_c$ denotes the space of 
continuously differentiable functions with compact support and
$\mathcal{N}$ is the set of eigenfunctions of the Stokes operator,
as defined in \cite{RobRodSad16}, Theorem 2.24.
(Caution: $\alpha_k \neq a_k$ in (\ref{eq:d}).)

Following \cite{RobRodSad16}, Def. 3.3, a weak solution of the Navier-Stokes
equation is defined as
\begin{enumerate}
\item
$u_0 \in \mathcal{H}$,

\item
$u \in L_\infty(0,T;\mathcal{H}) \cap L_2(0,T;\mathcal{H}_1) $,

\item
$u$ satisfies weak conservation of energy 
with $\phi \in \widetilde{{D_\sigma}}(\mathbb{T}^3)$ as test functions,
according to \cite{RobRodSad16}, Def. 3.1.
\begin{equation}
\label{eq:energy-weak}
\begin{split}
& \int_0^s - (u,\partial_t \phi )
+ \int_0^s \nu (\nabla u,\nabla \phi)
+ \int_0^s ((u \cdot \nabla u),\phi) \\
& = (u_0,\phi(0)) - (u(s),\phi(s)) \\
\end{split}
\end{equation}
for a.e. $s > 0$.
\end{enumerate}

Consequences of the definition, items 1--3,  of a weak solution
are given in  \cite{RobRodSad16}. 

Since $\widetilde{\mathcal{D}}_\sigma (\mathbb{T}^3) \subset 
L^1(0,t;\mathcal{V})$, we are free to use $L^1(0,T;\mathcal{V)}$
as the space of cutoffs. This fact motivates the definition of
$\mathcal{W}$ given earlier.

Using the reasoning of \cite{RobRodSad16} Exercise 3.4, $u$ solves the
Navier-Stokes equation
weakly, relative to test functions $\phi \in L^1(0,T;\mathcal{V}$).
Thus $u$ solves the Navier-Stokes equation in $\mathcal{W}$.

By \cite{RobRodSad16}, Lemmas 3.4 and 3.7,
\begin{equation}
\label{eq:sol-reg}
\begin{split}
& \partial_t u \in L^{4/3}(0,T);\mathcal{V}') \ , \\
& (u \cdot \nabla u) \in L^{6/5}(0,T;\mathcal{V}')
\subset L^{4/3}(0,T;\mathcal{V}') \ . \\
\end{split}
\end{equation}

\subsection{Second Order Phase Transitions}
\label{sec:sec-order}

Decay processes, referring to the entire turbulent process and
not to individual states within it, are subsumed in the generic notion
of a second order phase transition, whose ultimate rest state is the
hydrostatic state of a fluid at rest.

The second order phase transition occurs in four main steps, as indicated
here in a formal analysis and with the
first three steps studied numerically 
in \cite{HsuKauGli22}, based on the analysis of DNS simulation data.

The first step is the decay of enstrophy, from a highly stirred (high vorticity)
state to a state of near zero enstrophy. The second step is decay in the
energy mode, driven by the cascade of energy states to smaller wavelength,
where they ultimately vanish, and ending in the vacuum state of fully
developed turbulence. The third step is driven by viscous dissipation,
proceeding in a fourth step by radiation of thermal energy
to a regime in thermal equilibrium with an ambient temperature
and only thermally residual molecular motion.

The second order phase transitions for the Euler equation
have the first two and fourth steps, while lacking the 
viscosity related third step. 

The central point of this paper is that these transitions are achieved
discontinuously (with non-smoothness of the solution). Each non-smooth 
transition is realized by the removal of one or more 
degrees of freedom from the solution. 

With all these degrees of freedom removed, the solution is globally smooth
in time. Such smooth solutions are realized as the solution means of a
(general) weak Navier-Stokes solution.

\subsection{Fourier Analysis}
\label{sec:har}

Any $u \in \mathcal{H}$, the $L_2$ space of turbulent fields $u$,
has the Fourier expansion
\begin{equation}
\label{eq:energy-a}
u = \sum_{k \in \mathbb{Z}^3} \hat{u}_k e^{2\pi i x/L} \ ,
\end{equation}
with $L$ the length of one side of the periodic domain,
subject to the incompressibility condition
\begin{equation}
\label{eq:energy-t}
k \cdot \hat{u}_k = 0
\end{equation}
for all $k \in \mathbb{Z}^3$. Then
\begin{equation}
\label{eq:energy-c}
L^3  \sum_{ {k \in \mathbb{Z}^3}} |\hat{u}_k|^2 = |u|^2 < \infty \ .
\end{equation}

\paragraph*{Vector Fourier Analysis}

Vector Fourier analysis is described by a sum of representations
of the rotation group, $O(3)$. 

Vector Fourier analysis is achieved by adding a label $E$ or $\pm Z$
to the Fourier coefficients $k$ in (\ref{eq:energy-t}). The labels
define basis elements of the Lie algebra $\mathfrak{o}^3$ of the
global gauge group for a fluid velocity field.

\paragraph*{Young's Inequality}

Young's inequality \cite{FoiManRos04}, eq. (4.17), pg. 20.  states that
\begin{equation}
\label{eq:young}
ab \leq \frac{a^p}{p} + \frac{b^{p'}}{p'}
\end{equation}
for all $a,b > 0$ and for $p,p'$ Holder conjugate indices.

\paragraph*{Solid Harmonics}

Customary Fourier harmonics are
augmented with a magnitude $r$. The complete expansion has the basis
labeled $r,k$. On restriction to the unit sphere, $r$ takes on
two values, $r = \pm 1$. The plus sign selects the solid spherical 
harmonics solution that is regular at the origin. The minus sign
selects the solid spherical Fourier solution that is regular at infinity.

Due to restrictions to energy analysis subspaces, the solid Fourier
expansion is not used in this paper.

\subsection{Navier-Stokes Term by Term Regularity}
\label{sec:ns-reg}

If $u$ is a weak solution of the Navier-Stokes equation, then 
\begin{enumerate}
\item $u \in L^\infty(0,T;\mathcal{H}) \bigcap L^2(0,T;\mathcal{V})$,
\item $\partial_t u \in L^{4/3}(0,T;\mathcal{V}')$,
\item $(u \cdot \nabla )u \in L^{6/5}(0,T; \mathcal{V'})
\subset L^{4/3}(0,T;\mathcal{V}')$,
\item $u \in C(0,\infty;\mathcal{V}')$ with $C$
denoting the space of continuous functions,
\item $\partial_t u - \Delta u + (u \cdot \nabla) u = 0\in \mathcal{V}'$
as an element of $L^{4/3}(0,T;\mathcal{V}')$
\end{enumerate}
by \cite{RobRodSad16}, Lemma 2 (with following discussion) for item 1,
Lemma 3.7 for item 2, Lemma 3.4 for item 3 and Definition 3.1.
Item 4 is established in \cite{RobRodSad16}, pg. 76.
This is a consequence of item 3. Item 5 is a consequence of items 1 -- 3.

\begin{list}{}{\leftmargin=\parindent\rightmargin=0pt}
\item
\textbf{Proposition \ref{sec:ns-reg}.1 }
The terms in Item 5 are the terms of the Navier-Stokes equation. Each term is
individually spatially in $\mathcal{V}'$ and temporally in $L^{4/3}$. 
\end{list}
{\it Proof.}
As stated above.
\hspace*{\fill}$\square$

\subsection{The Schwartz Space}
\label{sec:sch}

The Schwartz space $\mathcal{S}$ of smooth vector fields is defined as
\begin{equation}
\label{eq:sch}
\mathcal{S} = C^\infty(0,T;C^\infty(\mathbb{T}^3))
\end{equation}
with its customary topology, and its dual space is
\begin{equation}
\label{eq:sch-dual}
\mathcal{S}' = C^{\infty'}(0,T;C^{\infty'}(\mathbb{T}^3))
\end{equation}

$\mathcal{S}$ and $\mathcal{S}'$ are closed under multiplication, so that
polynomials of elements of $\mathcal{S}$ are also in $\mathcal{S}$
and polynomial elements of $\mathcal{S}'$ are also in $\mathcal{S}'$.

\subsection{ Dissipation Operators}
\label{sec:diss}

The various forms of turbulent and viscous dissipation are summarized
by (\ref{eq:nu-NL}).

\begin{equation}
\label{eq:nu-NL}
\begin{split}
&\nu_{t,NL} = \langle \sum_j u_j \partial_j u_i \rangle \\
& \nu_{t,Temp} = \langle u_i,- \partial_t u_i \rangle \\
& \nu_t = \nu_{t,NL} + \nu_{t,Temp} \\ 
&\nu_{Tot} = \nu_{t,Temp}  +  \nu_{t,NL}  + \nu \|\nabla u \|^2 \ . \\
\end{split}
\end{equation}

The turbulent fluctuations without the expectation values are referred to
as the integrands:

\begin{equation}
\label{eq:nu-NL-int}
\begin{split}
&\nu_{t,NL}^\mathrm{int} = \sum_j u_j \partial_j u_i \\
& \nu_{t,Temp}^\mathrm{int} = u_i- \partial_t u_i \\
& \nu_t^\mathrm{int} = \nu_{t,NL}^\mathrm{int} + \nu_{t,Temp}^\mathrm{int} \ .\\ 
\end{split}
\end{equation}

The solid mechanics fixed space expansion variable $r = 0$
is regular at the origin.
Each of the individual terms in (\ref{eq:nu-NL}) is nonnegative for an
entropy production maximizing $u$.

\paragraph*{Non-Viscous Dissipation}

An entropy weak solution of the Navier-Stokes equation will dissipate
{\it non-negative energy} in a {\it non-viscous} manner.
This non-viscous dissipation
is a positive addition to its viscous dissipation. Referring to 
Fig.~\ref{fig:blowup}, the orange $\nu_t$ dissipation curve must always
be above the yellow viscous decay curve. This is point 3 of the
three contradictory facts of Sec.~\ref{sec:reason}.

In estimating polynomials of the dissipation operators (\ref{eq:nu-NL}),
cross terms are bounded by Young's inequality, so that bounds of polynomials
in the individual terms are sufficient.

\paragraph*{Natural Function Spaces}

The proofs of analyticity and blowup occur in $\mathcal{V}^*$.
After multiplication by a test function $\phi$ to allow for weak convergence
and after entropy production maximization, the proofs are based in
$\mathcal{V}^*$ and $\mathcal{W}$.
Proofs are simplified in comparison to $\mathcal{H}$ based proofs.
The proof of finite time blowup depends on the $\mathcal{V}^*$ time universal
decay of the turbulent fluctuations (\ref{eq:nu-NL-int}) $O(t^{-5/6})$ 
and $O(t^{-4/3})$.  These decay rates are strictly larger than the
viscous decay $O(t^{-1/2})$
and result in a finite time blowup if the initial data contains
turbulent fluctuations.

The evaluation of the large time limit involves solution means, with the
analysis occuring in the space $\mathcal{H}$. The mean is shown to have zero
entropy. It is SRI and its Navier-Stokes evolution equation reduces to the
heat equation. The mean initial data is also SRI and in view of
this restriction, is a proper subset of all allowed initial data.

\subsection{$\mathbb{P}$ in Vector Fourier Harmonics}
\label{sec:p-sph}

The construction of an entropy production maximizing weak solution of the
Navier-Stokes equation is too singular to address directly. This result 
is needed for the proof of non-negativity of the $\nu_t$ and analyticity.

The present section is an analysis of the Leray projector 
$\mathbb{P}$ as an operator in the space of vector Fourier harmonics.

\paragraph*{Navier-Stokes as Gauge Field}

The gauge group for the Navier-Stokes equation 
is the group of global $O(3)$ symmetries. Any global
$O(3)$ transformation leaves all its statistical properties unchanged.

Here we refer to the formalism of 
\cite{IvePhi08}, eq. (12) and (13). 

\paragraph*{$\mathbb{P}$ in Vector Fourier Harmonics}

The \cite{IvePhi08} formalism is useful
as it includes a definition of the divergence, and thus of  the Leray
projector $\mathbb{P}$. The needed formulas are \cite{IvePhi08}  eq. (23) and (24). From these,
we see that $\mathbb{P}$ is defined within a single $k$ scalar Fourier
harmonics space and in this space it is a combination of 
its components in the three vector Fourier harmonics spaces.

\subsection{Initial Data}
\label{sec:II}

The function spaces $L^{4/3}(0,T;\mathcal{V}^*$) and $\mathcal{W}$ describe
weak solutions of the Navier-Stokes equation. It is natural to work
within these spaces. Specificly, the spatial function space $\mathcal{V}^*$
and the temporal space $L^{4/3}$ play key roles. The 
space $\mathcal{H}$ is not used in the blowup analysis.

\begin{list}{}{\leftmargin=\parindent\rightmargin=0pt}
\item
\textbf{Proposition \ref{sec:II}.1 }
A weak solution of the Navier-Stokes equation is in $C(0,T,\mathcal{V^*})$.
\end{list}
{\it Proof.}
The use of $\mathbb{P}$ in the individual Fourier harmonics 
subspaces was justified above.
\hspace*{\fill}$\square$

\paragraph*{No Global Drift}

\begin{list}{}{\leftmargin=\parindent\rightmargin=0pt}
\item
\textbf{Proposition~\ref{sec:II}.2}
Let $u$ be a weak solution of the Navier-Stokes equation with smooth
initial data and a stirring function $f$. Then there is a solution
$\tilde{u}$ and a stirring function $\tilde{f}$ both with zero global
average. $\tilde{u}$ and $\tilde{f}$ are constructed from $u$ and $f$
by a spatially constant solution of the Navier-Stokes equation.
Thus we assume that $u$ and $f$ have a zero spatial average.
\end{list}
{\it Proof.} This is a simple computation.
See also \cite{FoiManRos04}, Pg. 31 and  \cite{RobRodSad16} exercise 3.8.
\hspace*{\fill}$\square$

In the remainder of the paper, the Navier-Stokes solution is assumed to
have a zero spatial average.

\paragraph*{SRI and NSRI}

A summation over repeated tensor indices is denoted SRI. A sum over 
tensor indices in which the SRI term is omitted is denoted NSRI.

Important for this paper are the NSRI $\nu^\mathrm{int}_t$
turbulent fluctuation terms of (\ref{eq:nu-NL-int}).

\begin{list}{}{\leftmargin=\parindent\rightmargin=0pt}
\item
\textbf{Proposition \ref{sec:II}.4 }
(i) For $u$ a weak solution of the Navier-Stokes equation,
\begin{equation}
\label{eq:u2-transport}
\begin{split}
&\int_{t_1}^{t_2} \langle \nabla u,\nabla \phi \rangle
+ \int_{t_1}^{t_2} \langle (u \nabla) u, \phi \rangle  \\
& = \langle u(t_1) , \phi \rangle - \langle u(t_2),\phi\rangle \\
\end{split}
\end{equation}
for all cutoff functions $\phi$ and a.e. $t_2 \geq t_1$ and a.e. $t_1 \geq 0$
including $t_1 = 0$.

(ii) $\partial_ t u \in L^{4/3}(0,T;\mathcal{V}^*)$

(iii) One can modify $u$ on a set of measure zero so that 
$u$ is a continuous function of $t$ with values in
$\mathcal{V}^*$ 
and $u$ is still a solution of its transport equation.
\end{list}
{\it Proof.}
The proof is based on \cite{RobRodSad16}, Lemmas 3.6, 3.7.
\hspace*{\fill}$\square$

The viscous dissipation term in the above propositions
is necessarily SRI due to restrictions
on the possible forms of viscous dissipation. Thus there is no subscript $i$
in the viscous dissipation term.

\section{The Entropy Principle}
\label{sec:ent}

\paragraph*{Entropy and Enstrophy}

As a vector field, the entropy of the Navier-Stokes solution
is maximized for each basis
element of the three dimensional gauge group Lie algebra $\mathfrak{o}^3$. 
In simpler terms, this maximizes conventional (scalar) entropy and
enstrophy. The entropy is defined as the log volume of a surface of
constant energy plus the log volume of a surface of constant enstrophy within
a surface of constant energy.

The maximum production of entropy is a fundamental law of physics.
We refer to this set of ideas as the entropy principle, for short.

The entropy principle as used in this paper states the existence of entropy
production maximizing turbulent fluctuation integrands $\nu_t^\mathrm{int}$
and solutions $u$ in the energy sector.

For a quantitative comparison, it is also necessary to specify the function
space in which the maximum occurs. 

The function space is $\mathcal{W}$ with common
initial conditions. 

The entropy principle $\nu_t^\mathrm{int}$ have an entropy production
which is maximal relative to any $\nu_t^\mathrm{int}$
in the same function space and initial conditions. It is likewise 
maximal for the solution $u$.

This principle is used with cutoffs only. Here cutoffs
are applied to the $\nu_t^\mathrm{int}$
turbulent fluctuation terms and to the solution $u$.
As such, the principle is
finite dimensional, and its existence has been established by multiple authors.
Thus the key issue is the function space in which the convergence occurs.

The existence of an entropy and enstrophy production
maximizing state is the main result of \cite{CheGliSai24b}.
See also \cite{GliLazChe20a}.
However, this existence theorem is not directly usable here due to the
function space in which the maximum entropy occurs.

\paragraph*{Entropy and Laws of Physics}

It is necessary to specify the physical laws assumed for the current
analysis. Assuming isothermal, single fluid physics, entropy is time
independent, while the energy satisfies a heat equation with decay of energy.

Other physical assumptions, such as thermal physics and compressibility
allow time dependent entropy.

For the isothermal physics considered here, the entropy has a dependence on the
constant temperature of the solution. A convenient normalization is to
set the constant temperature to infinity, defining the microcanonical ensemble.

\paragraph*{Weak Energy Conservation}

A weak form of the (exact) conservation of energy was given by M. Lee
(private communication) with an eight line proof.
The proof assumes a sum over repeated indices (SRI). We extend this result 
to turbulent fluctuation dissipation and to the NSRI case
(which is defined as no SRI). The extended conservation law
plays a key role in the paper.

Let $u$ be a weak solution of the Navier-Stokes equation. From
\cite{RobRodSad16}, Lemma 3.6,
\begin{equation}
\label{eq:weak-en}
\partial_t u -\nu \Delta u + (u\cdot\nabla)u  = 0 \in \mathcal{V'} 
\end{equation}
for a.e. $t$; $u \in \mathcal{V'}$ for a.e. $t$.

The weak exact conservation of energy is formulated as 

\begin{list}{}{\leftmargin=\parindent\rightmargin=0pt}
\item
\textbf{Lemma \ref{sec:ent}.1 }
The energy is conserved under weak limits, up to viscous losses.
\end{list}
{\it Proof.}
Take the dot product of 
(\ref{eq:weak-en}) with the solution $u$ and vacuum expectation values.
The pairing $u \in \mathcal{V}$, $\partial_t u \in \mathcal{V}'$ yields
\begin{equation}
\label{eq:weak-en-1}
\langle \partial_t u, u \rangle + \langle (u \cdot \nabla)u,u \rangle
= -\nu \| (\nabla u)^2 \| 
\end{equation}
and
\begin{equation}
\label{eq:cons-en}
\partial_t \langle u,u \rangle + \langle u,-\partial_t u \rangle 
 + \langle (u \cdot \nabla)u,u \rangle
= -\nu\|(\nabla u)^2\|    
\end{equation}
for a.e. $t$.

The cubic term vanishes by \cite{RobRodSad16}, Lemma 3.2, assuming SRI.
The inertial term $\nu_{t,Temp}$ is also zero assuming SRI, again by
\cite{RobRodSad16}, Lemma 3.7 and the following discussion. 
\hspace*{\fill}$\square$

Lemma~\ref{sec:ent}.1 is limited to the energy analysis as defined
by the Navier-Stokes global gauge Lie algebra $\mathfrak{o}^3$.
There is no mention of conservation of enstrophy $\pm Z$.
Thus all that follows in this paper is restricted to the energy
analysis. 

The use of Lemma~\ref{sec:ent}.1 for $u$ and $\nu_t$ requires additional
entropy production maximization.
As discussed in \cite{RobRodSad16}, an application without modification
of the $v_{NL}$ proof encounters the space $V^*$. This space is
defined in \cite{RobRodSad16}, Chapt. 1 as the Banach space dual to $V$.

\begin{list}{}{\leftmargin=\parindent\rightmargin=0pt}
\item
\textbf{Theorem \ref{sec:ent}.1 }
(a) An entropy principle solution of the Navier-Stokes equation exists in
$\mathcal{W}$
which maximizes the entropy production of $u$ and of $\nu_t^\mathrm{int}$
relative to the entropy production of $u$ and of $\nu_t^\mathrm{int}$ 
defined by any weak solution of the Navier-Stokes equation in the same
function space and with the same initial conditions.

(b) The entropy principle $u$ and $\nu_t^\mathrm{int}$ are elements of 
$\mathcal{W}$, with decay $O(e^{-4\nu t/3})$. 

(c) The entropy principle $\nu_t^\mathrm{int}$
is non-negative and at least equal to the decay of the solution $u$.

The norm of entropy principle $u$ in $\mathcal{W}$
and its norm in $\mathcal{W}$ is smaller than
the norm of $\nu_t^\mathrm{int}$ in $\mathcal{W}$.
\end{list}
{\it Proof.}
Weak solutions with a nonzero spatial average are excluded.

The proof is based on cutoffs, with convergence 
in $\mathcal{W}$ uniform in the cutoffs.

The cutoff problem is finite dimensional and the existence of 
an entropy production maximizing solutions for the cutoff $\nu_t$
and $ u $ is classical.

The fluctuation is a real number, which can be positive,
zero or negative. The entropy principle excludes the negative choice,
so that the entropy principle fluctuation is non-negative.
This proves (a).

The entropy principle is applied to the set $\mathcal{W}$.
Since $L^{4/3}(0,T;\mathcal{V}^* L^1(0,T;\mathcal{V}))$
was already identified as a
limit point, the entropy principle selection for $u$ within the set 
$\mathcal{V}^*L^1(0,T;\mathcal{V})$ must dissipate turbulent
fluctuations at least  according to 
$L^{4/3}(0,T;\mathcal{V}^*L^1(0,T;\mathcal{V}))$. 
This proves (b).

The norm of $u \in \mathcal{W}$
is smaller than the norm of $\nu_t^\mathrm{int} \in \mathcal{W}$.
This statement is derived from the proof of Lemma~\ref{sec:ent}.1.
All terms of Lemma~\ref{sec:ent}.1 are SRI. The entropy principle
$u$ is the sum of the non-negative $\nu_t^\mathrm{int}$ and a positive
viscous decay. Thus the entropy principle $u$ is strictly smaller than
$\nu_t^\mathrm{int}$ if $\nu_t^\mathrm{int}$ is not zero.
This proves (c)
\hspace*{\fill}$\square$

\paragraph*{Legendre Transforms}

The energy dissipation principle asserts the existence of a weak solution
of the Navier-Stokes equation which dissipates at least as much 
$\nu_t^\mathrm{int}$ or $u$ as any comparison solution. 

The detailed formulation of this principle requires specification of
the component of the solution being compared (the $\nu_t^\mathrm{int}$
and $u$) and the function space ($\mathcal{W}$) with common initial
conditions.

\begin{list}{}{\leftmargin=\parindent\rightmargin=0pt}
\item
\textbf{Theorem \ref{sec:ent}.2 }
A weak solution of the Navier-Stokes equation exists 
in the space $L^{4/3}(0,T;\mathcal{V}^*L^1(0,T;\mathcal{V}))$
which maximizes the dissipation rate 
of the fluctuation $\nu_t^\mathrm{int}$ and the solution $u$
relative to that of the fluctuation the $\nu_t^\mathrm{int}$
or the solution $u$
defined by any weak solution of the Navier-Stokes equation in
$\mathcal{W}$ with the same initial data.
 \end{list}
{\it Proof.}
Entropy is related to energy and enstrophy by a double Legendre transform,
first in energy, then in enstrophy and then repeated for three basis
elements $T^a$ of the global gauge Lie algebra $\mathfrak{o}^3$,
or in simpler terms,
 repeated for the entropy and enstrophy variables.

The Legendre transforms of the entropy (in the thermodynamic space of entropy,
energy and enstrophy) depends on convexity of the entropy, first as
a function of energy, and then at fixed energy, as a function of the
enstrophy.

The convexity proofs start with the convexity of the logarithm function,
first of the energy and then of the enstrophy.

As entropy is defined as
the log volume of a surface of constant energy, convexity of the
entropy as a function of the energy follows. At fixed energy, the log enstrophy
is convex as a function of the enstrophy. Entropy, in this case, is defined
as the log volume of a surface of constant enstrophy (and constant energy).
\hspace*{\fill}$\square$

\paragraph*{One-Sided Limits}

One-sided time limits of turbulent fluctuating quantities $\nu_t^\mathrm{int}$
are a regularity property of the time evolution.
There are distinct reasons for continuity from the left (limits from
earlier time) and continuity from the right (limits from later time).

It is known that the solution is strong (in $\mathcal{V}$) except for 
an isolated and possibly empty set of singular times $t_j$ at which times a
blowup occurs.

Approaching $t_j$ from the left (earlier times),
the turbulent fluctuation $\nu_t^\mathrm{int}L^1(0,T;\mathcal{V})$
is analytic with
$\nu_t^\mathrm{int}  \in \mathcal{W}$ decreasing in time
and always above the viscous decay $O(e^{-\nu t/2})$. Thus the left
(from earlier times) is limit is convergent. The same conclusions are
reached in \cite{FoiManRos04}. In addition, from the same reference,
is the existence of a countable (possibly empty) sequence of singular
times $t_j$ at which the solution  fails to be smooth. 

Approaching $t_j$ from the right (later times), the solution is strong,
in $\mathcal{V}$, and continuous.

\section{Analyticity}
\label{sec:III}

The complexified Navier-Stokes equation in $\mathcal{V}^*$ is
\begin{equation}
\label{eq:complex-ns}
\frac{du}{d \zeta} + (-\nu\Delta) u + B(u)_\mathbb{C} = 0 \in \mathcal{V}^* \ 
\end{equation}
\cite{FoiManRos04}, eq. (8.1).
Here we write $(-\Delta)(u_1 + iu_2) = (-\Delta) u_1 + i(-\Delta) u_2$.

The complexified bilinear form $B_\mathbb{C}$  is written as sum of its
real and imaginary parts.
The same splitting of the
other complex terms in (\ref{eq:complex-ns}) into real and imaginary terms 
is used, with individual estimates for the real and imaginary terms.

Due to the periodic boundary conditions, the Stokes operator has been
reduced to $-\nu\Delta$. This simple identity reduces estimates of the
complex term to estimates of the corresponding real terms, as in
\cite{FoiManRos04}. Analyticity for the Stokes operator in $\mathcal{H}$,
established in \cite{FoiManRos04} extends to analyticity in $\mathcal{V}^*$.

Analyticity of $u$ in $\mathcal{V}^*$ is elementary to prove directly
(as a bound on the nonlinear term) and extends to analyticity for 
$\nu_t^\mathrm{int}$.

For example see Sec.~\ref{sec:ent}
for time independent bounds for $\nu_t^\mathrm{int}$ in $\mathcal{V}^*$.

The proof of analyticity for the turbulent fluctuation integrands
is based on multiplication by a test function from $L^1(0,T;V^*)$,
the time universal turbulent fluctuation dissipation of Sec.~\ref{sec:ent},
with an entropy principle hypothesis, and non-negativity of the 
$\nu_t^\mathrm{int}$.

The result is UV finite norms in $\mathcal{W}$ for $\nu_{t,NL}^\mathrm{int}$
and $\nu_{t,Temp}^\mathrm{int}$
with $\mathcal{W}$ continuity uniform in $t$ in a limit as the cutoffs are
removed, from Theorem~\ref{sec:ent}.1.

\begin{list}{}{\leftmargin=\parindent\rightmargin=0pt}
\item
\textbf{Proposition \ref{sec:III}.1 }
The real domain bounds of Proposition~\ref{sec:ent}.1 extend to the
complex domain. 
\end{list}
{\it Proof.}
The real domain analysis is extended to the
complex domain. As in \cite{FoiManRos04}, the complex domain analysis
is repeated for the real and imaginary parts and the two bounds added.
The result is still finite in $\mathcal{W}$, uniformly in $t$, with a
larger but still finite norm.

For the turbulent dissipation terms, the same reduction of the proof
to the real case only applies. 

The proof is completed by the Theorem~\ref{sec:ent}.1 bounds in
$\mathcal{W}$
on the turbulent dissipation terms that are uniformly nonzero in
$t$ as the cutoffs are removed.
\hspace*{\fill}$\square$

\paragraph*{Turbulent Initial Data}

\begin{list}{}{\leftmargin=\parindent\rightmargin=0pt}
\item
\textbf{Theorem \ref{sec:III}.1 }
Let $u$ be a weak solution of the Navier-Stokes equation with initial
data orthogonal with zero spatial average. For example, let $u$ be a LH solution
with zero spatial average.

$u$ is a strong solution except for possible isolated
singular times $t_j \neq 0$.
\end{list}
{\it Proof.}
According to (\ref{eq:sol-reg}), the initial conditions are
achieved continuously in $\mathcal{V}^*$. 
$u \in \mathcal{V}$ is a strong solution.
\hspace*{\fill}$\square$

\section{Blowup}
\label{sec:blowup}

\begin{list}{}{\leftmargin=\parindent\rightmargin=0pt}
\item
\textbf{Definition~\ref{sec:blowup}.1}
Initial data $u_0$ is called laminar (non-turbulent) if 
it is SRI
Such basis functions are proportional to a mean value.

The data is otherwise defined to be turbulent.
\end{list}

The numerical program of \cite{Hou23} documents a rapid 
near total blowup in enstrophy, followed by a slower blowup in
the energy, to complete the blowup process. 

The initial data is chosen with
turbulent initial conditions.
The LH weak solution defined by this data lies in $\mathcal{V}^*$.
We prove blowup for this weak solution. To stay within the energy analysis
of Lemma~\ref{sec:ent}.1, the
turbulent initial data are energy variables. Enstrophy variables are avoided.

\begin{list}{}{\leftmargin=\parindent\rightmargin=0pt}
\item
\textbf{Theorem~\ref{sec:blowup}.1} An entropy maximizing solution of the
Navier-Stokes equation with turbulent energy mode initial conditions 
has only a finite interval of regularity. It is not globally smooth in time.
\end{list}

{\it Proof.}
Assume the solution is within an interval of regularity starting from
its initial conditions. 

Initial conditions are chosen with the $E$ mode strictly non-zero,
non-SRI. In other words, they are chosen with non-zero NSRI energy modes.
These initial conditions are smooth since they lie 
within the span of a finite number of eigenvalues of the Stokes operator.

Consider the $E$ analysis of the subspace of connected Feynman diagrams of
length $j = 3$. The vector indices each have dimension 3. There is an
index for the component of the fluid vector velocity and a component for
the spatial derivative of this velocity, as in the formula for the
strain rate. The total number of indices is 9. The SRI subspace has
dimension 3. The remaining 6 dimensions are NSRI, and include $\nu_t$.

Blowup occurs in the case of non-zero initial data in this $E$ 
NSRI subspace.

$\nu_{t}(s) \in \mathcal{W}$ decreases as a function of $s$ at least at
the rate $O(e^{-4\nu t/3})$
and is always greater than the solution $u$ decay curve.

The viscous decay rate for $u$ is $O(e^{-\nu t/2})$.

For a finite value of $s$, $\nu_t(s)$ 
is smaller than the viscous decaying $u$,
with the decay rate $O(e^{-\nu t/2})$, since $4/3 > 1/2$.

Thus at this time, $\nu_t$ falls
below the smooth solution viscous decay allowed minimum value for $u$.
An endpoint to
the epoch of regularity has been reached. This time is $T^* < \infty$.
$T^*$ is the time of blowup.

The discontinuity (and resulting non-smoothness of the solution) 
occurs in the blowup mode, which is an NSRI
expansion coefficient for some connected Feynman diagram.

Some connected Feynman diagram is discontinuous at $t = t_j$,
which demonstrates the lack of smoothness.

The length $j$ of the diagram has $j = 1$ excluded as this is the
conserved solution mode. $j = 2$ is allowed, and defines $\nu_v^\mathrm{int}$.
The discontinuity could also occur in higher $j$ diagrams.
\hspace*{\fill}$\square$

This set of inequalities is illustrated in Fig.~\ref{fig:blowup}.

The solution $u$ is continuous and not the source of the discontinuity.

A discontinuity occurs in some connected Feynman diagram, as it is
converted from turbulent to viscous.
A finite unit of turbulent fluctuation has transitioned through a
range of small length scales, at or below the Kolmogorov scale,
where viscosity is the primary dissipative mechanism, to become totally
viscous dissipation. Thus the turbulent dissipation is discontinuous. 
An NSRI mode has become a SRI mode, discontinuously. 
By an entropy principle, the discontinuity is complete, meaning that
the entire NSRI mode has become SRI.

\section{The Mean}
\label{sec:mean-min}

The mean $\mathcal{M}(u)$, defined in the vector Fourier
space $\mathcal{H}$ is a sum of three separate means.
One is defined in the energy expansion and the other two are the $\pm$ oriented
enstrophy expansions. The $\pm$ vorticity orientation reduces to two separate
formulas for the mean $\mathcal{M}(u)$. 

\subsection{Entropy Production Minimization}
\label{sec:ent-min}

We note that the total
entropy is equal to the sum of two individual entropies,
defined by the $r$ (fixed space) and the fixed time expansion variables.

A Legendre transform relates fixed time energy dissipation to the 
fixed time entropy. 

All proofs are of the form of uniform properties of cutoff quantities.
The Legendre transforms are
established in the cutoff case. Thus we can reason in terms of 
energy dissipation. 

\paragraph*{The Minimum Entropy Principle}

The minimum entropy principle (when valid for a specific solution)
states that a nonnegative polynomial function of the
solution has a minimum rate of
entropy production relative to the same polynomial function
of any other weak solution of the Navier-Stokes equation with the same
initial conditions. 

\paragraph*{A Key Distinction}

The key distinction 
between the two (opposite) cases of entropy production maximizing
vs. minimizing is the fact that for an
entropy production minimizing solution, the turbulent dissipation
rate $\nu_t$ is zero. 

The blowup in Fig. 1 
changes its character: the diagonal sloping line
generated by the time derivative is now horizontal, and never intersects
the horizontal line below it defined by the energy.
Thus there is no end to the duration of the strong solution, and the
mean solution is smooth for all space and time.

\begin{list}{}{\leftmargin=\parindent\rightmargin=0pt}
\item
\textbf{Remark \ref{sec:ent-min}.1 }
The minimum entropy principle is used in this paper to study Navier-Stokes
solutions and also
nonnegative polynomial functions of $\nu_{t,NL}$ and $\nu_{t,Temp}$.

The fixed time energy 
is intrinsically nonnegative. A nonnegative entropy principle polynomial
must have nonnegative even coefficients and zero odd coefficients. 
\end{list}

\begin{list}{}{\leftmargin=\parindent\rightmargin=0pt}
\item
\textbf{Theorem \ref{sec:ent-min}.1 }
Any element of $\mathcal{H}$ is expressed uniquely in terms of
its fixed time and fixed space expansion variables.

The mean $\mathcal{M}(u)$ is the sum (normalized by division by the number of
terms in the sum) of 
\begin{enumerate}
\item the $r$ configuration space (magnitude equal $1$) mean,
\item the loop expansion phase space mean defined at fixed spatial coordinates,
\item the mean defined at fixed $t$.
\end{enumerate}
The mean has entropy zero.
\end{list}

The proof is placed in Sec.~\ref{sec:app}.

\subsection{Regularity of the Mean}
\label{sec:mean}

For a nonnegative polynomial function $p(u)$ of a weak solution of the
Navier-Stokes
equation, its average value $\mathcal{M}(|p(u)|)$ is defined as the ratio 
$\langle |p(u)| \rangle/ \langle 1 \rangle$. In simple terms it is the
sum of values divided by the number of terms in the sum. 

Here we study $\mathcal{M}(\nu_{t,Temp})$ and $\mathcal{M}(\nu_{t,NL})$.
Products of distributions are
not a concern, as we introduce cutoffs (so that  all distributions become
smooth functions which can be multiplied) and then prove results that are
uniform as the cutoffs are removed.

The complete mean is the sum of the means defined by the fixed time and 
fixed space means.

\begin{list}{}{\leftmargin=\parindent\rightmargin=0pt}
\item
\textbf{Theorem \ref{sec:mean}.1 }
\begin{equation}
\label{eq:mean-3}
\mathcal{M}(\nu_{t,Temp}) = 0 \in L^{4/3}(0,T;\mathcal{V}') \subset \mathcal{S}'
\end{equation}
\begin{equation}
\label{eq:mean-2}
\mathcal{M}(\nu_{t,NL}) = 0 \in L^{4/3}(0,T;\mathcal{V}') \subset \mathcal{S}'
\end{equation}
\end{list}

\begin{list}{}{\leftmargin=\parindent\rightmargin=0pt}
\item
\textbf{Theorem \ref{sec:mean}.2 }
\begin{equation}
\label{eq:mean-5}
\nu_{t,Temp}(\mathcal{M}(\nu_{t,Temp})) = 0 \in \mathcal{S}'
\end{equation}
\begin{equation}
\label{eq:mean-4}
\nu_{t,NL}(\mathcal{M}(\nu_{t,NL})) = 0 \in L^{4/3}(0,T;\mathcal{H}_{-3})
\subset{S}' \ .
\end{equation}
\end{list}

Intuitively, the evaluation is zero because the turbulent fluctuation of 
a turbulent fluctuation is zero.

Proofs are placed in Sec.~\ref{sec:app}.

\section{The Mean is a Smooth Solution of the Navier-Stokes Equation}
\label{sec:2}

The main result of this section is the
construction of smooth solutions of the Navier-Stokes equation with smooth
non-turbulent initial data. The smooth solution is
the mean of a weak solution. Its entropy is zero, so that it can also be
characterized as an entropy production minimizing solution. 
We show that

\begin{enumerate}

\item Initial conditions: The initial conditions of the mean of a
weak solution of the Navier-Stokes
equation with initial data $u^0 \in \mathcal{H}(\mathbb{T}^3)$
(the space of $L_2$ velocity fields defined on the cube $\mathbb{T}^3$
with periodic boundary conditions) conditions are
assumed in a continuous fashion, in the sense that 
\begin{equation}
\label{eq:ic}
\mathcal{M}(u^0) = \lim_{t \rightarrow 0} \mathcal{M}(u(t)) \ .
\end{equation}
in the topology of $\mathcal{H}(\mathbb{T}^3)$.

All vector fields are assumed to be divergence-free. 

\item 
The mean of a weak solution of the Navier-Stokes equation is a strong solution.
It is smooth for all time. 
\end{enumerate}

\begin{list}{}{\leftmargin=\parindent\rightmargin=0pt}
\item
\textbf{Theorem \ref{sec:2}.1 }
The mean is smooth for all positive times.

The mean assumes its initial conditions (\ref{eq:ic}) in $\mathcal{H}$.

The mean of a weak solution $u$ of the Navier-Stokes equation in $\mathcal{H}$
is a solution of the Navier-Stokes equation in $\mathcal{H}$.
\end{list}
{\it Proof.}
By Theorem~\ref{sec:mean}.1, 
$\mathcal{M}(\nu_{t,NL}(u)) = 0 \in \mathcal{H}_{-2} \subset \mathcal{S}'$.
and
$\mathcal{M}(\nu_{t,Temp}(u)) = 0 \in \mathcal{H}_{-2} \subset \mathcal{S}'$.

To find the equation which the mean satisfies, we reason
directly. $u$ satisfies the equation
\begin{equation}
\label{eq:ns-weak}
\partial_t u + \nu_{t,Temp} + \nu_{t,NL}(u) = \nu \Delta u \in \mathcal{S}'\ .
\end{equation}

Since the mean $\mathcal{M}$ is linear,
it is applied to each of the terms of (\ref{eq:ns-weak}) individually, in
the space $\mathcal{S}'$.

$\mathcal{M}(\nu_{t,NL}(u))$ and $\mathcal{M}(\nu_{t,Temp}(u))$ are each 
the 0 element of $\mathcal{S}'$.

By Theorem~\ref{sec:mean}.1, there remains the equation
\begin{equation}
\label{eq:ns-mean}
\partial_t \mathcal{M}(u) = \nu \Delta \mathcal{M}(u) \in \mathcal{S}' \ .
\end{equation}

This is the equation is also the heat equation. Continuity of initial
conditions in $\mathcal{H}$ is based on an assumed the global spatial
translation invariance of the solution. This property is preserved by
taking a mean.

It is
well known that the solutions in $\mathcal{S}'$ of the heat equation with 
$\mathcal{S}'$ initial conditions are smooth for positive times, belonging to 
$\mathcal{S}$. In view of the $\mathcal{H}$ continuity of the data, the
initial data is in $\mathcal{H}$. The first two assertions are proven.

To prove the third assertion, recall the 0 evaluation in $\mathcal{S}'$ of
$\nu_{t,NL}(\nu_{t,NL}(\mathcal{M}{)u})$ 
and of $\nu_{t,Temp}(\mathcal{M}(\nu_{t,Temp}))$ according to
Theorem~\ref{sec:mean}.2. 

Thus the cutoff version of (\ref{eq:ns-mean}) is the
cutoff Navier-Stokes equation in $\mathcal{S}'$.

As inferred above, $\mathcal{M}(u) \in \mathcal{S}$ is smooth for positive
times, so that 
(\ref{eq:ns-mean}) is the Navier-Stokes equation in $\mathcal{S}$.
\hspace*{\fill}$\square$.

\section{Appendix}
\label{sec:app}

\subsection{Theorem~\ref{sec:ent-min}.1}
\label{sec:app1}

\begin{list}{}{\leftmargin=\parindent\rightmargin=0pt}
\item
\textbf{Theorem \ref{sec:app1}.1 }
The minimum entropy solution of the Navier-Stokes equation is unique.
It is an element of $\mathcal{V}'$.
\end{list}
{\it Proof:}
By items 2,3, $\nu_{t,Temp}\in\mathcal{V}'$ and
$\nu_{t,NL} \in \mathcal{V}'$.

By general theories of positive definite forms, uniqueness of a minimum is
a result of strict convexity.

Strict convexity of the energy dissipation is a result of the
viscous dissipation term $\nu \Delta u$, which sets a strictly positive
lower bound on the dissipation as $\nu \Delta u$, with a gradient in
$\mathcal{V}'$. This completes the proof for a solution in $\mathcal{V}'$.
\hspace*{\fill}$\square$

{\it Proof of Theorem~\ref{sec:ent-min}.1}
A mean is not a random variable. It is the expectation value of a random 
variable.

The first statement is elementary.

The mean of any of the summands 1 -- 4  in the statement of the theorem
is defined as
\begin{equation}
\label{eq:mean-energy}
\frac
{\sum \mathrm{summand}}
{\mathrm{number}\; \mathrm{of}\; \mathrm{terms}} \ .
\end{equation}
In the cutoff cases where this formula is used, the sum contains only a 
finite number of terms.

The $r$ phase space entropy production and $r$ phase space energy dissipation
are evaluated at fixed spatial
coordinates while the fixed time entropy production and energy dissipation are
evaluated with the time coordinate fixed. 

The $r$ configuration space dissipation is a one dimensional.
Uniqueness of this (one dimensional)
minimum dissipation state follows.

The $r$ conjugate momentum space is the expansion space of the conjugate
momenta $\pi$, called the loop expansion. Uniqueness of the entropy minimizing
loop expansion variable is a consequence of Theorem~\ref{sec:app1}.1.

The volume of the entropy minimizing solution is its $L_1$ expectation
value, also equal to its $L_1$ norm. The volume of
the solution divided by its $L_1$ norm is 1. With the unity temperature
normalization, the entropy is the $L_1$ norm of this ratio.

By the definition of the
entropy, the $r$ entropy of the mean is $\log 1 = 0$.

Next we specialize to a fixed time.

The denominator in (\ref{eq:mean-energy}) is the
number of energy dissipation terms in the numerator.

All steps as in the $r$ case can be followed. 
The energy dissipation is intrinsicly nonnegative, to reduce the
solution space to one dimension. To arrive at a one dimensional space,
we use Theorem~\ref{sec:ent-min}.2, which states that the fixed time entropy
minimum is unique, or
equivalently, that the fixed time minimum energy dissipation state is unique.

(Related analysis focused on entropy production maximization has been 
presented in joint work with H. Said and G.-Q. Chen.)

The conclusion is then the same as in the $r$ case.
The entropy of the fixed time mean $\mathcal{M}(u)$ is zero. 

The plus sign is chosen since the vacuum expectation value 
is nonnegative. The rest of the proof follows as in
the $r$ case.

Thus
\begin{equation}
\label{eq:unit-vol}
\langle M(u) \rangle =
\langle \mathcal{M}(u) \rangle_\mathcal{H} = 1 \ .
\end{equation}

Since the volume of $\langle \mathcal{M}u) \rangle_\mathcal{H} = 1$,
the entropy (the log of this volume) of $\mathcal{M} = 0$.

Entropy is always 0 or larger.
Thus the mean is entropy minimizing.
\hspace*{\fill}$\square$

\begin{list}{}{\leftmargin=\parindent\rightmargin=0pt}
\item
\textbf{Corollary \ref{sec:app1}.1 }
The mean satisfies the minimum entropy principle.
\end{list}
{\it Proof.}
Since $\mathcal{S}$ and $\mathcal{S}'$ are closed under multiplication,
the polynomial powers needed here are automatic. Thus the proof follows 
as in Theorem~\ref{sec:ent-min}.1.
\hspace*{\fill}$\square$

\subsection{Proof of Theorem~\ref{sec:mean}.1}

By items 2,3, $\nu_{t,Temp}\in\mathcal{V}^*$ and
$\nu_{t,NL} \in \mathcal{V}^*$.
The mean of $\nu_{t,Temp}$ is at least as regular as $\nu_{t,Temp}$.
Likewise, the mean of $\nu_{t,NL}$ is at least as regular as $\nu_{t,NL}$.

To justify the value 0 for the mean, consider item 2 of Sec.~\ref{sec:ns-reg}
for (\ref{eq:mean-3}) and item 3 for (\ref{eq:mean-2}).

The time derivative $\partial_t u$ occurs exclusively in the
$r$ expansion.  

The nonlinear term $\langle u_i, - \partial_t u_i \rangle$ occurs only in the
loop expansion. The term is evaluated at fixed spatial
coordinates. As a result of Proposition~\ref{sec:ent}.1, 
the energy in the $r$ loop expansion is conserved up to the
amount dissipated by the viscous term $\Delta u$. 

Think of the mean as a sum of terms. Each term is a weak solution,
to which Lemma~\ref{sec:ent}.1 is applicable. 
preserved by taking a mean. Thus the $r$ loop expansion terms contributing to 
$\mathcal{M}(\nu_{t,temp})$ is a sum of zeros in $\mathcal{V}^*$.

The $r$ term (\ref{eq:mean-3}) selects loop expansion
which is regular at the origin.

The proof of (\ref{eq:mean-3}), based on the $r$ configuration space and the
$r$ loop expansion is complete.

Next, we prove (\ref{eq:mean-2}).

The nonlinear term $(u \cdot \nabla)u$ occurs exclusively
in the fixed time expansion.
At fixed time, the time derivative and
the $r$ expansion are removed from item 3 of Sec.~\ref{sec:ns-reg}.

The mean is a sum of terms to
which Lemma~\ref{sec:ent}.1 is applied. For each summand, the result is
zero. 
The energy expansion is conserved exactly
up to the amount dissipated by the viscous term $\Delta u$. 
The sum of zeros is still a zero.
The nonlinear terms contributing to
$\mathcal{M}(\nu_{t,NL})$ sum to zero. 

This term is
nonnegative, which was to be demonstrated. The quadratic polynomial
$\nu_{t,NL}$ has a nonnegative quadratic coefficient, to
complete the proof of (\ref{eq:mean-2}).
\hspace*{\fill}$\square$

\subsection{Proof of Theorem~\ref{sec:mean}.2}

To prove (\ref{eq:mean-5}) analyticly, consider
\begin{equation}
\label{eq:Tmean-Tmean}
\begin{split}
& \nu_{t,Temp}(\mathcal{M}(\nu_{t,Temp})) 
= \nu_{t,Temp}\mathcal{M}[\langle (u_i, - \partial_t u_i) \rangle] \\
& = \Big \langle 
\mathcal{M}[\langle (u_i, - \partial_t u_i) \rangle],  
- \partial_t \mathcal{M}[\langle (u_i, - \partial_t u_i) \rangle],
\Big \rangle \\
= 0 \ .  \\
\end{split}
\end{equation}

Here the value is 0 because it is the fluctuation of a mean, while
the mean is a number and not a statistical quantity; its fluctuations are
zero.

There is a problem with the direct determination of the function space
that this zero is located in. But introducing cutoffs, everything is
regular, so that the cutoff values of (\ref{eq:Tmean-Tmean}) are zero in 
$\mathcal{H}$ for example.

Next we prove (\ref{eq:mean-4}) analytically.  Consider
\begin{equation}
\begin{split}
\label{eq:mean-mean}
&\nu_{t,NL}\mathcal{M}(\nu_{t,NL}) =
\langle (u \cdot \nabla) \cdot \nabla u\mathcal{M}\langle (u\cdot \nabla u
\rangle\rangle  \\
& = \langle ((u \cdot \nabla) \cdot \nabla u)\langle (u \cdot\nabla)(u\cdot\nabla)u \rangle \rangle   \ . \\
\end{split}
\end{equation}
The value is again zero and for the same reason. The mean is a number
and does not fluctuate. Its nonlinear fluctuations are zero.
The zero for the cutoff version of (\ref{eq:mean-mean}) is in $\mathcal{H}$.
\hspace*{\fill}$\square$

\section{Acknowledgements and Disclaimer}
We thank G.-Q. Chen, H. Said, M. Lee, T. Wallstrom, T. Spencer
and A. Rahimyar for helpful comments.

Point72 disclaimer:

The information, views, and opinions expressed herein are solely Jarret
Petrillo's and do not necessarily represent the views of Point72 or its
affiliates.  Point72 and its affiliates are not responsible for, and did not
verify for accuracy, any of the information contained herein.

\bibliographystyle{siam}
\bibliography{refs.bib}
\end{document}